\documentclass{amsart}
\usepackage{amssymb}
\usepackage{graphicx}

\newtheorem{thm}{Theorem}
\newtheorem{cor}{Corollary}
\newtheorem{lem}{Lemma}

\newtheorem{rem}{Remark}

\newtheorem{conj}{Conjecture}
\theoremstyle{definition}
\newtheorem{example}[equation]{Example}
\newtheorem{prob}[equation]{Problem}

\newcommand{\A}{{\mathcal A}}

\newcommand{\ID}{{\mathbb D}}

\newcommand{\D}{{\mathbb D}}

\newcommand{\real}{{\operatorname{Re}\,}}

\def\be{\begin{equation}}
\def\ee{\end{equation}}

\newcommand{\bee}{\begin{enumerate}}
\newcommand{\eee}{\end{enumerate}}

\newcommand{\blem}{\begin{lem}}
\newcommand{\elem}{\end{lem}}
\newcommand{\bthm}{\begin{thm}}
\newcommand{\ethm}{\end{thm}}
\newcommand{\bcor}{\begin{cor}}
\newcommand{\ecor}{\end{cor}}
\newcommand{\beg}{\begin{example}}
\newcommand{\eeg}{\end{example}}
\newcommand{\begs}{\begin{examples}}
\newcommand{\eegs}{\end{examples}}
\newcommand{\bdefe}{\begin{defin}}
\newcommand{\edefe}{\end{defin}}
\newcommand{\bprob}{\begin{prob}}
\newcommand{\eprob}{\end{prob}}
\newcommand{\bei}{\begin{itemize}}
\newcommand{\eei}{\end{itemize}}

\newcommand{\bcon}{\begin{conj}}
\newcommand{\econ}{\end{conj}}
\newcommand{\bcons}{\begin{conjs}}
\newcommand{\econs}{\end{conjs}}
\newcommand{\bprop}{\begin{propo}}
\newcommand{\eprop}{\end{propo}}
\newcommand{\br}{\begin{rem}}
\newcommand{\er}{\end{rem}}
\newcommand{\brs}{\begin{rems}}
\newcommand{\ers}{\end{rems}}
\newcommand{\bo}{\begin{obser}}
\newcommand{\eo}{\end{obser}}
\newcommand{\bos}{\begin{obsers}}
\newcommand{\eos}{\end{obsers}}
\newcommand{\bpf}{\begin{pf}}
\newcommand{\epf}{\end{pf}}
\newcommand{\ba}{\begin{array}}
\newcommand{\ea}{\end{array}}
\newcommand{\beq}{\begin{eqnarray}}
\newcommand{\beqq}{\begin{eqnarray*}}
\newcommand{\eeq}{\end{eqnarray}}
\newcommand{\eeqq}{\end{eqnarray*}}

\begin{document}

\title[  Toeplitz and symmetric Toeplitz determinants ]{ Toeplitz and symmetric Toeplitz determinants for inverse functions of certain classes of univalent
functions }

\author[M. Obradovi\'{c}]{Milutin Obradovi\'{c}}
\address{Department of Mathematics,
Faculty of Civil Engineering, University of Belgrade,
Bulevar Kralja Aleksandra 73, 11000, Belgrade, Serbia.}
\email{obrad@grf.bg.ac.rs}

\author[N. Tuneski]{Nikola Tuneski}
\address{Department of Mathematics and Informatics, Faculty of Mechanical Engineering, Ss. Cyril and
Methodius
University in Skopje, Karpo\v{s} II b.b., 1000 Skopje, Republic of North Macedonia.}
\email{nikola.tuneski@mf.ukim.edu.mk}



\begin{abstract}
In this paper we investigate  Toeplitz and symmetric Toeplitz determinants of inverse functions for some classes of univalent functions and improve some previous results.
\end{abstract}

\maketitle

\section{Introduction and preliminaries}

Let ${\mathcal A}$ be the class of functions that are analytic in the unit disk $\D:= \{ |z| < 1 \}$
and  are normalized such that
$f(0)=0= f'(0)-1$, i.e.,
\begin{equation}\label{eq-1}
f(z)=z+a_2z^2+a_3z^3+\cdots.
\end{equation}
By ${\mathcal S}$ we denote the class of functions from ${\mathcal A}$ which are univalent in $\D$. Properties of univalent functions are studied on the general class $\mathcal{S}$ or on its subclasses. Some subclasses if the class of univalent functions are the class of functions with bounded turning, the class of starlike  functions and the class of convex functions, defined respectively in the following way
\[
\begin{split}
\mathcal{R}&=\left\{f\in\mathcal{A}: {\rm Re}f'(z)>0,\, z\in\ID \right\},\\
\mathcal{S}^{\star}&=\left\{f\in\mathcal{A}: {\rm Re}\frac{zf'(z)}{f(z)}>0,\, z\in\ID \right\},\\
\mathcal{C}&=\left\{f\in\mathcal{A}: {\rm Re}\left[1+\frac{zf''(z)}{f'(z)}\right]>0,\, z\in\ID \right\}.
\end{split}
\]

\medskip

In the recent time, estimation of the modulus of the Toepliz and the symmetric Toepliz determinant over the general class $\mathcal{S}$ or its subclasses attracts significant interest of the researchers working in the field. For a function $f\in \A$ of the form \eqref{eq-1} and positive integers $q$ and $n$, the  Toepliz determinant is defined by
\[
        T_{q,n}(f) = \left|
        \begin{array}{cccc}
        a_{n} & a_{n+1}& \ldots& a_{n+q-1}\\
        \overline{a}_{n+1}&a_{n}& \ldots& a_{n+q-2}\\
        \vdots&\vdots&~&\vdots \\
        \overline{a}_{n+q-1}& \overline{a}_{n+q-2}&\ldots&a_{n}\\
        \end{array}
        \right|,
\]
where $\overline{a}_k = \overline{a_k}$, and as given in \cite{firoz} (see \cite{MONT-2025}), and the symmetric Toepliz determinant is defined by
\[
        T^{s}_{q,n}(f) = \left|
        \begin{array}{cccc}
        a_{n} & a_{n+1}& \ldots& a_{n+q-1}\\
        a_{n+1}&a_{n}& \ldots& a_{n+q-2}\\
        \vdots&\vdots&~&\vdots \\
        a_{n+q-1}& a_{n+q-2}&\ldots&a_{n}\\
        \end{array}
        \right|.
\]

\medskip

Thus, the second  Toepliz determinant is
\begin{equation}\label{eq-2}
T_{2,1}(f)= 1-|a_2|^2
\end{equation}
and the third is
\begin{equation}\label{eq-3}
T_{3,1}(f) =  \left |
        \begin{array}{ccc}
        1 & a_2& a_3\\
        \overline{a}_2 & 1& a_2\\
        \overline{a}_3 & \overline{a}_2& 1\\
        \end{array}
        \right | = 2\, \real (a_2^2\overline{a}_3) -2|a_2|^2-|a_3|^2+1 .
\end{equation}
The concept of  Toeplitz matrices plays an important role in functional analysis, applied mathematics as well as in physics and technical sciences (for more details see \cite{ye}).

\medskip

In their paper \cite{MONT-2021} the authors gave the next results for the class S:

\medskip

\noindent
\textbf{Theorem A.}
\emph{If $f\in\mathcal{S}$, then
\[-3\le T_{2,1}(f)\le 1 \quad\mbox{and}\quad -1\le T_{3,1}(f)\le 8. \]
All inequalities are sharp.}

\medskip

From the previous definition we easily obtain that
\begin{equation}\label{eq-4}
\begin{split}
T^{s}_{2,2}(f)=&a_{2}^{2}-a_{3}^{2},\\
T^{s}_{2,3}(f)=&a_{3}^{2}-a_{4}^{2},\\
T^{s}_{3,1}(f)=&1-2a_{2}^{2}+2a_{2}^{2}a_{3}-a_{3}^{2},\\
T^{s}_{3,2}(f)=&(a_{2}-a_{4})(a_{2}^{2}-2a_{3}^{2}+a_{2}a_{4}).
\end{split}
\end{equation}

\medskip

In \cite{firoz} the author proved the following estimate of some symmetric Toepliz determinants over the general class of univalent functions.

\medskip

\noindent
\textbf{Theorem B.}
\emph{If $f\in\mathcal{S}$, then 
\[|T^{s}_{2,2}(f)|\leq 13 ,\quad |T^{s}_{2,3}(f)|\leq 25 ,\quad |T^{s}_{3,1}(f)|\leq 24. \]
All inequalities are sharp.}

\medskip

In this paper our main focus will be on finding estimates of the modulus of the  Toeplitz and symmetric Toeplitz determinants for inverse functions of the functions in the class $\mathcal{S}$ and its subclasses $\mathcal{R}$, $\mathcal{S}^{\star}$ and $\mathcal{C}$.

\medskip

Namely, for every univalent function in $\D$, exists its inverse at least on the disk with radius 1/4 (due to
the famous Koebe's 1/4 theorem). If the inverse has an expansion
\begin{equation}\label{eq-5}
f^{-1}(w) = w+A_2w^2+A_3w^3+\cdots,
\end{equation}
then, by using the identity $f(f^{-1}(w))=w$, from \eqref{eq-1} and \eqref{eq-5} we receive

\begin{equation}\label{eq-6}
\begin{split}
A_{2}&=-a_{2}, \\
A _{3}&=-a_{3}+2a_{2}^{2} , \\
A_{4}&= -a_{4}+5a_{2}a_{3}-5a_{2}^{3}.
\end{split}
\end{equation}

In the proofs we will use the following result from \cite{PS-1981}.

\begin{lem}\label{lem-1}
Let $\Omega$ denote the class of analytic functions $\omega: \D\rightarrow \D$ of the form
\begin{equation}\label{eq-7}
\omega(z)=c_{1}z+c_{2}z^{2}+\cdots.
\end{equation}
If $\omega\in \Omega$ is given by \eqref{eq-7}, then
$$|c_1|\le 1,\quad |c_2|\le1-|c_1|^2 $$
and
$$ \left|c_3+ \mu c_1 c_2 +\nu c^{3}_1\right|\leq|\nu|,$$
where $\mu$ and $\nu$ are real and $(\mu,\nu)$ belongs to $D_{4},\,D_{6},$ or $D_{7}$, where
\[
\begin{split}
D_{4}&=\left\{(\mu,\nu): |\mu|\geq \frac{1}{2},\, \nu\leq -\frac{2}{3}(|\mu|+1)\right\}, \\
D_{6}&=\left\{(\mu,\nu): |\mu|\geq 4,\, \nu\geq \frac{2}{3}(|\mu|-1)\right\}, \\
D_{7}&=\left\{(\mu,\nu): 2\leq|\mu|\leq 4,\, \nu\geq\frac{1}{12}(\mu^{2}+8)\right\}. 
\end{split}
\]
\end{lem}

\medskip

\section{Toeplitz determinants for inverse functions}

The next theorem is valid for Toeplitz determinants for inverse functions
of univalent functions.

\bthm\label{25-th-1}
If $f\in \mathcal{S}$, then
$$T_{2,1}(f^{-1})=T_{2,1}(f)\quad \mbox{and}\quad T_{3,1}(f^{-1})=T_{2,1}(f).$$
\ethm

\begin{proof}
From the relations  \eqref{eq-2}, \eqref{eq-3} and \eqref{eq-6}, we have

$$T_{2,1}(f^{-1})=1-|A_{2}^{2}|=1-|-a_{2}^{2}|=T_{2,1}(f),$$
and
\[
\begin{split}
T_{3,1}(f^{-1})=&2\, \real (A_2^2\overline{A}_3) -2|A_2|^2-|A_3|^2+1 \\
=&2\, \real \left[(-a_2)^2(-\overline{a}_3+2(\overline{a}_2)^{2})\right] -2|-a_2|^2-\left|-a_{3}+2a_{2}^{2}\right|^2+1\\
=& 2\, \real \left(-a_2^2\overline{a}_3+2|a_{2}|^{4}\right)-2|a_{2}|^{2}- \left[|a_{3}|^{2}+4|a_{2}|^{4}-4\real(a_2^2\overline{a}_3)\right] +1\\
=&2\, \real (a_2^2\overline{a}_3) -2|a_2|^2-|a_3|^2+1 \\
=&T_{3,1}(f).
\end{split}
\]
\end{proof}

From Theorem A and Theorem \ref{25-th-1} we have the next
\begin{cor}
 If $f\in\mathcal{S}$, then
\[-3\le T_{2,1}(f^{-1})\le 1 \quad\mbox{and}\quad -1\le T_{3,1}(f^{-1})\le 8. \]
All inequalities are sharp with extremal functions that reach the boundaries of the estimates being the Koebe function $k(z)=\frac{z}{(1-z)^{2}}$ and its inverse
$k^{-1}(w)=w-2w^{2}+5w^{3}-14w^{4}+\cdots $, and the function $f_{1}(z)=\frac{z}{1-z+z^{2}}=z+z^{2}-z^{4}+\cdots $ with its inverse $f^{-1}_{1}(w)=w-w^{2}+2w^{3}-4w^{4}+\cdots$.
\end {cor}

\medskip

We note that Theorem \ref{25-th-1} implies the appropriate results for Toeplitz determinants of inverse functions in case for other subclasses of univalent functions (see \cite{cudna,Kowalczyk-20}).

\medskip

\section{Symmetric Toeplitz determinants for inverse functions}

For the class $\mathcal{S}$ we have the next statement.
\bthm\label{25-th-2}
If $f\in\mathcal{S}$, then 
\[|T^{s}_{2,2}(f^{-1})|\leq 29 ,\quad |T^{s}_{2,3}(f^{-1})|\leq 221 \quad \mbox{and}\quad |T^{s}_{3,1}(f^{-1})|\leq 24. \]
All inequalities are sharp.
\ethm

\begin{proof}
Since for $f\in\mathcal{S}$ we have $|a_{2}|\leq2$, $|a_{3}|\leq3$, $|a_{4}|\leq4 $ and $|a_{3}-a_{2}^{2}|\leq1$
(see \cite{duren}), using \eqref{eq-6}, we easily obtain 
\[
\begin{split}
|A_{2}| &=|-a_{2}|\leq2, \\
|A_{3}| &=|-a_{3}+2a_{2}^{2}|\leq|a_{2}|^{2}+ |a_{3}-a_{2}^{2}|\leq5,\\
|A_{4}| &= \left|-a_{4}+5a_{2}(a_{3}-a_{2}^{2})\right|\leq|a_{4}|+ 5|a_{2}||a_{3}-a_{2}^{2}|\leq14.
\end{split}
\]
We note that all the previous estimates are sharp as the inverse of the Koebe function
$k^{-1}(w)=w-2w^{2}+5w^{3}-14w^{4}+\cdots $ shows.

\medskip

Now, using these estimates for $|A_{2}|$, $|A_{3}|$, $|A_{4}|$, and \eqref{eq-4}, we have
\[
\begin{split}
|T^{s}_{2,2}(f^{-1})| &=|A_{2}^{2}-A_{3}^{2}|\leq |A_{2}|^{2}+|A_{3}|^{2}\leq 29,\\
|T^{s}_{2,3}(f^{-1})| &=|A_{3}^{2}-A_{4}^{2}|\leq |A_{3}|^{2}+|A_{4}|^{2}\leq 221.
\end{split}
\]
Also,
\[
\begin{split}
T^{s}_{3,1}(f^{-1})=& 1-2A_{2}^{2}+2A_{2}^{2}A_{3}-A_{3}^{2}\\
=& 1-2(-a_{2})^{2}+2(-a_{2})^{2}(-a_{3}+2a_{2}^{2})-(-a_{3}+2a_{2}^{2})^{2}\\
=&1-2a_{2}^{2}+2a_{2}^{2}a_{3}-a_{3}^{2}\\
=& T^{s}_{3,1}(f),
\end{split}
\]
which implies $|T^{s}_{3,1}(f^{-1})|=|T^{s}_{3,1}(f)|\leq 24,$ where we used Theorem B.

\medskip

The estimates of the modulus of the symmetric Toeplitz determinants above are the best possible as the inverse function
 $f^{-1}_{2}(w)= w-2iw^{2}-5w^{3}+14iw^{4}+\cdots $   of the function
$f_{2}(z)=\frac{z}{(1-iz)^{2}}$ shows.
\end{proof}

In their paper \cite{Hadi-2025}, the authors gave the next three theorems.

\medskip

\noindent
\textbf{Theorem C.}
\emph{If $f\in \mathcal{R}$, then
\begin{itemize}
\item[(a)]  $|T^{s}_{2,2}(f^{-1})|\leq 7.22 $;
\item[(b)]  $|T^{s}_{2,3}(f^{-1})|\leq 168.694$;
\item[(c)]  $|T^{s}_{3,1}(f^{-1})|\leq3.88 $;
\item[(d)]  $|T^{s}_{3,2}(f^{-1})|\leq 64.79$.
\end{itemize}
These estimates are sharp (except inequality $T^{s}_{3,2}$) for when
\[ f(z)=\frac{1+iz}{1-iz}=1+2iz-2z^{2}-2iz^{3}+2z^{4}+\cdots  .\]}

\noindent
\textbf{Theorem D.}
\emph{If $f\in \mathcal{S}^{\star}$, then
\begin{itemize}
\item[(a)]  $|T^{s}_{2,2}(f^{-1})|\leq 51 $;
\item[(b)]  $|T^{s}_{2,3}(f^{-1})|\leq 116.33$;
\item[(c)]  $|T^{s}_{3,1}(f^{-1})|\leq24 $;
\item[(d)]  $|T^{s}_{3,2}(f^{-1})|\leq 650.56$.
\end{itemize}
All sharp estimates are derived from the inverse of
\[ f(z)=\frac{z}{(1-iz)^{2}}=z+2iz^{2}-3z^{3}-4iz^{4}+\cdots  .\]}

\noindent
\textbf{Theorem E.}
\emph{If $f\in \mathcal{C}$, then
\begin{itemize}
\item[(a)]  $|T^{s}_{2,2}(f^{-1})|\leq 2.7 $;
\item[(b)]  $|T^{s}_{2,3}(f^{-1})|\leq 10.27$;
\item[(c)]  $|T^{s}_{3,1}(f^{-1})|\leq4 $;
\item[(d)]  $|T^{s}_{3,2}(f^{-1})|\leq 7.24$.
\end{itemize}
All sharp estimates are derived from the inverse of $ f(z)=\frac{z}{1-iz} .$}

\medskip

We want to stress that all results given in Theorem C, Theorem D and Theorem E are not correct apart from the cases (c) in all these theorems.

\medskip

Using different approach, in our next theorems the correct results will be given.

\bthm\label{25-th-3}
If $f\in \mathcal{R}$, then we have
\begin{itemize}
\item[($i$)]  $|T^{s}_{2,2}(f^{-1})|\leq\frac{25}{9}=2.77\ldots  $ ;
\item[($ii$)]  $|T^{s}_{2,3}(f^{-1})|\leq \frac{233}{36}=6.472\ldots $;
\item[($iii$)]  $|T^{s}_{3,2}(f^{-1})|\leq \frac{817}{108}=7.5448\ldots  .$
\end{itemize}
All estimates are sharp for the function defined by
$$ f'(z)=\frac{1+iz}{1-iz}=1+2iz-2z^{2}-2iz^{3}+2z^{4}+\cdots ,$$
i.e., for the function
$$ f(z)=z+iz^{2}-\frac{2}{3}z^{3}-\frac{1}{2}iz^{4}+\cdots  .$$
Its inverse function is
$$f^{-1}(w)=w-iw^{2}-\frac{4}{3}w^{3}+\frac{13}{6}iw^{4}+\cdots  .$$
\ethm

\begin{proof}
Since $f\in \mathcal{R}$ is equivalent to
$$f'(z)=\frac{1+\omega(z)}{1-\omega(z)} (=1+2\omega(z)+2\omega^{2}(z)+\cdots ,$$
where $\omega \in \Omega.$
Using the notations for $f$ and $\omega$ given by \eqref{eq-1} and \eqref{eq-7}, and equating the coefficients in previous relation, we receive
\be\label{eq-8}
\begin{split}
  a_2=& c_1, \\
  a_3=& \frac{2}{3}(c_2+c_1^2), \\
  a_4=& \frac{1}{2}(c_3+2c_1c_2+c_1^3).
\end{split}
\ee
From the relations \eqref{eq-6} and \eqref{eq-8}, after some simple calculations we have
\begin{equation}\label{eq-9}
\begin{split}
A_{2}&=-c_{1}, \\
A_{3}&=-\frac{2}{3}c_{2}+\frac{4}{3}c_{1}^{2} , \\
A_{4}&=-\frac{1}{2}\left(c_{3}-\frac{14}{3}c_{1}c_{2}+\frac{13}{3}c_{1}^{3}\right).
\end{split}
\end{equation}
Now, by applying Lemma \ref{lem-1}, we get $|A_{2}|\leq1$,\,
$$|A_{3}|\leq \frac{2}{3}|c_{2}|+\frac{4}{3}|c_{1}|^{2}
\leq\frac{2}{3}(1-|c_{1}|^{2})+\frac{4}{3}|c_{1}|^{2}\leq \frac{4}{3},$$
while
$$|A_{4}|= \frac{1}{2}\left|c_{3}-\frac{14}{3}c_{1}c_{2}+\frac{13}{3}c_{1}^{3}\right|
\leq\frac{1}{2}\cdot \frac{13}{3}=\frac{13}{6}.$$
In last inequality we used $(\mu,\nu) = (-14/3, 13/3)\in D_{6}$. So,
\begin{equation}\label{eq-10}
|A_{2}|\leq1,\quad |A_{3}|\leq\frac{4}{3}\quad \mbox{and}\quad |A_{4}|\leq \frac{13}{6}.
\end{equation}

\medskip

It is evident that all these estimates are sharp as the function $f$ given in the statement of this theorem shows.

\medskip

For proving the estimates (i), (ii), and (iii), from  \eqref{eq-4} and \eqref{eq-10} we easily obtain
$$|T^{s}_{2,2}(f^{-1})|=|A_{2}^{2}-A_{3}^{2}|\leq |A_{2}|^{2}+|A_{3}|^{2}\leq\frac{25}{9}$$
and
$$|T^{s}_{2,3}(f^{-1})|\leq |A_{3}|^{2}+|A_{4}|^{2}\leq\frac{233}{36}.$$
Also, since
\[
\begin{split}
\left|A_{2}^{2}-2A_{3}^{2}+A_{2}A_{4}\right| 
=& \left|c_{1}^{2}-\frac{8}{9}c_{2}^{2}
+\frac{1}{2}c_{1}\left(c_{3}+\frac{22}{9}c_{1}c_{2}-\frac{25}{9}c_{1}^{3}\right)\right| \\
\leq& |c_{1}|^{2}+\frac{8}{9}|c_{2}|^{2}
+\frac{1}{2}|c_{1}|\left|c_{3}+\frac{22}{9}c_{1}c_{2}-\frac{25}{9}c_{1}^{3}\right| \\
\leq& |c_{1}|^{2}+\frac{8}{9}(1-|c_{1}|^{2})^{2}+\frac{1}{2}\cdot \frac{25}{9}\\
=&\frac{8}{9}-\frac{7}{9}|c_{1}|^{2}+\frac{8}{9}|c_{1}|^{4}+\frac{25}{18}\\
\leq& 1+ \frac{25}{18}=\frac{43}{18}.
\end{split}
\]
In the third step above we applied Lemma \ref{lem-1} with $(\mu,\nu) = (22/9. -25/9)\in D_{4}$. 
Then finally
$$\left|T^{s}_{3,2}(f^{-1})\right| \leq \left(|A_{2}|+|A_{4}|\right) \cdot \left|A_{2}^{2}-2A_{3}^{2}+A_{2}A_{4}\right|\leq\left(1+\frac{13}{6}\right)\frac{43}{18}=\frac{817}{108}.$$

\medskip

All the results of this theorem are the best possible as the function given in its statement shows.
\end{proof}

\medskip

Now we will study the starlike functions.
    
\bthm\label{25-th-4}
If $f\in \mathcal{S}^{\star}$, then we have
\begin{itemize}
\item[($i$)]  $|T^{s}_{2,2}(f^{-1})|\leq29$;
\item[($ii$)]  $|T^{s}_{2,3}(f^{-1})|\leq 221$;
\item[($iii$)]  $|T^{s}_{3,2}(f^{-1})|\leq 416.$
\end{itemize}
All results are sharp.
\ethm
\begin{proof}
From the definition of the class $\mathcal{S}^{\star}$  we have that there exists a
a function $\omega\in \Omega$ such that
$$\frac{zf'(z)}{f(z)} = \frac{1+\omega(z)}{1-\omega(z)},$$
and from here
$$z f'(z)=\left[1+2\left(\omega(z)+\omega^{2}(z)+\cdots\right)\right]\cdot f(z).$$
From the last relation, as in the two previous theorems, by comparing the coefficients we receive
\be\label{eq-11}
\begin{split}
a_{2}&=2c_{1}, \\
a_{3}&=c_{2}+3c_{1}^{2},\\
a_{4}&=\frac{2}{3}\left(c_{3}+5 c_{1}c_{2}+6 c_{1}^{3}\right).
\end{split}
\ee
From relations \eqref{eq-6} and \eqref{eq-11} we obtain
\begin{equation}\label{eq-12}
\begin{split}
A_{2}&=-2c_{1}, \\
A_{3}&=-c_{2}+5c_{1}^{2} , \\
A_{4}&=-\frac{2}{3}\left(c_{3}-\frac{10}{3}c_{1}c_{2}+21c_{1}^{3}\right).
\end{split}
\end{equation}
Now, using Lemma \ref{lem-1}, we easily get
\begin{equation}\label{eq-13}
|A_{2}|\leq2,\quad |A_{3}|\leq 5\quad \mbox{and}\quad |A_{4}|\leq 14.
\end{equation}
We note that for the estimate of $|A_{4}|$ we used Lemma \ref{lem-1} with $(\mu,\nu)= (-10/3,21)\in D_{7}$).

\medskip

All the results in \eqref{eq-13} are sharp as the function $f_{2}$ given in Theorem \ref{25-th-2} shows.

\medskip

Further, using \eqref{eq-4} and \eqref{eq-13} we easily obtain
$$|T^{s}_{2,2}(f^{-1})|\leq |A_{2}|^{2}+|A_{3}|^{2}\leq29,$$
$$|T^{s}_{2,3}(f^{-1})|\leq |A_{3}|^{2}+|A_{4}|^{2}\leq 221,$$
and
\[
\begin{split}
|A_{2}^{2}-2A_{3}^{2}+A_{2}A_{4}|
=&\left|4c_{1}^{2}-2c_{2}^{2}
+\frac{4}{3}c_{1}\left(c_{3}+\frac{35}{3}c_{1}c_{2}-\frac{33}{2}c_{1}^{3}\right)\right| \\
\leq& 4|c_{1}|^{2}+2|c_{2}|^{2}
+\frac{4}{3}|c_{1}|\left|c_{3}+\frac{35}{3}c_{1}c_{2}-\frac{33}{2}c_{1}^{3}\right| \\
=&2+2|c_{1}|^{4}+22 \\
\leq& 26 .
\end{split}
\]
Here we used Lemma \ref{lem-1} with  $(\mu,\nu) = (35/3, -33/2) \in D_{4}$)).
So,
$$|T^{s}_{3,2}(f^{-1})| \leq \left(|A_{2}|+|A_{4}|\right) \cdot \left|A_{2}^{2}-2A_{3}^{2}+A_{2}A_{4}\right| \leq(2+14)26=416.$$

\medskip

The results are the best possible as the function $f(z)=\frac{z}{(1-iz)^{2}}$ and its  inverse function
 $f^{-1}(w)= w-2iw^{2}-5w^{3}+14iw^{4}+\cdots $ show.
\end{proof}

\medskip

At the end we cover the case of convex functions.

\bthm\label{25-th-5}
If $f\in \mathcal{C}$, then we have
\begin{itemize}
\item[($i$)]  $|T^{s}_{2,2}(f^{-1})|\leq2$;
\item[($ii$)]  $|T^{s}_{2,3}(f^{-1})|\leq 2$;
\item[($iii$)]  $|T^{s}_{3,2}(f^{-1})|\leq 4.$
\end{itemize}
All results are sharp.
\ethm
\begin{proof}

We apply the same method as in the two previous case. Namely, from the definition of the class $\mathcal{C}$  we have

$$1+\frac{zf''(z)}{f'(z)} = \frac{1+\omega(z)}{1-\omega(z)},$$
where $\omega\in\Omega$, and from here
$$\left(z f'(z)\right)'=\left[1+2\left(\omega(z)+\omega^{2}(z)+\cdots\right)\right]\cdot f'(z).$$

Using the notations \eqref{eq-1} and \eqref{eq-7}, and comparing the coefficients  in the previous relation, after some simple calculations, we obtain
\be\label{eq-14}
\begin{split}
a_{2}&=c_{1}, \\
a_{3}&=\frac{1}{3}\left(c_{2}+3c_{1}^{2}\right),\\
a_{4}&=\frac{1}{6}\left(c_{3}+5 c_{1}c_{2}+6 c_{1}^{3}\right).
\end{split}
\ee
From relations \eqref{eq-6} and \eqref{eq-14} we can get the coefficients of inverse functions:
\begin{equation}\label{eq-15}
\begin{split}
A_{2}&=-c_{1}, \\
A_{3}&=-\frac{1}{3}(c_{2}-3c_{1}^{2}) , \\
A_{4}&=-\frac{1}{6}(c_{3}-5c_{1}c_{2}+6c_{1}^{3}).
\end{split}
\end{equation}
By applying Lemma \ref{lem-1} we easily obtain
\begin{equation}\label{eq-16}
|A_{2}|\leq1,\quad |A_{3}|\leq 1\quad \mbox{and}\quad |A_{4}|\leq 1.
\end{equation}
From \eqref{eq-14} we get
\[
\begin{split}
\left|A_{2}^{2}-2A_{3}^{2}+A_{2}A_{4}\right|
=& \left| c_{1}^{2}-\frac{2}{9}c_{2}^{2}
+\frac{1}{6}c_{1}\left(c_{3}+3c_{1}c_{2}-6c_{1}^{3}\right)\right| \\
\leq& |c_{1}|^{2}+\frac{2}{9}|c_{2}|^{2}
+\frac{1}{6}|c_{1}|\left|c_{3}+3c_{1}c_{2}-6c_{1}^{3}\right| \\
\leq&|c_{1}|^{2}+\frac{2}{9}(1-|c_{1}|^{2})^{2}+\frac{1}{6}\cdot1\cdot6 \\
\leq& 2 .
\end{split}
\]
As in the previous theorem  we finally obtain the statement of the theorem. 

\medskip

The results are the best possible as the function
$f(z)=\frac{z}{1-iz}$ and its  inverse function
 $f^{-1}(w)= w-iw^{2}-w^{3}+iw^{4}+\cdots $ show.
\end{proof}

\end{document}